\documentclass[11pt]{amsart}
\usepackage{amssymb,mathrsfs}
\title{Holomorphic line bundles on a special abelian surface}

\begin{document}

\author{Jae-Hyun Yang}

\address{Yang Institute for Advanced Study
\newline\indent
Seoul 07989, Korea
\vskip 2mm
and
\vskip 2mm
Department of Mathematics
\newline\indent
Inha University
\newline\indent
Incheon 22212, Korea}

\email{jhyang@inha.ac.kr}

\newtheorem{theorem}{Theorem}[section]
\newtheorem{lemma}{Lemma}[section]
\newtheorem{proposition}{Proposition}[section]
\newtheorem{remark}{Remark}[section]
\newtheorem{definition}{Definition}[section]

\renewcommand{\theequation}{\thesection.\arabic{equation}}
\renewcommand{\thetheorem}{\thesection.\arabic{theorem}}
\renewcommand{\thelemma}{\thesection.\arabic{lemma}}
\newcommand{\BR}{\mathbb R}
\newcommand{\BQ}{\mathbb Q}
\newcommand{\BF}{\mathbb F}
\newcommand{\BT}{\mathbb T}
\newcommand{\BM}{\mathbb M}
\newcommand{\bn}{\bf n}
\def\charf {\mbox{{\text 1}\kern-.24em {\text l}}}
\newcommand{\BC}{\mathbb C}
\newcommand{\BZ}{\mathbb Z}


\begin{abstract}
We consider a special abelian surface $A_\Omega$ deduced from the work of
Tianze Wang, Tianqin Wang and Hongwen Lu \cite{WWL}. We study holomorphic line bundles over
a special abelian surface explicitly.
\end{abstract}

\maketitle

\newcommand\tr{\triangleright}
\newcommand\al{\alpha}
\newcommand\be{\beta}
\newcommand\g{\gamma}
\newcommand\gh{\Cal G^J}
\newcommand\G{\Gamma}
\newcommand\de{\delta}
\newcommand\e{\epsilon}
\newcommand\z{\zeta}
\newcommand\vth{\vartheta}
\newcommand\vp{\varphi}
\newcommand\om{\omega}
\newcommand\p{\pi}
\newcommand\la{\lambda}
\newcommand\lb{\lbrace}
\newcommand\lk{\lbrack}
\newcommand\rb{\rbrace}
\newcommand\rk{\rbrack}
\newcommand\s{\sigma}
\newcommand\w{\wedge}
\newcommand\fgj{{\frak g}^J}
\newcommand\lrt{\longrightarrow}
\newcommand\lmt{\longmapsto}
\newcommand\lmk{(\lambda,\mu,\kappa)}
\newcommand\Om{\Omega}
\newcommand\ka{\kappa}
\newcommand\ba{\backslash}
\newcommand\ph{\phi}
\newcommand\M{{\Cal M}}
\newcommand\bA{\bold A}
\newcommand\bH{\bold H}
\newcommand\D{\Delta}

\newcommand\Hom{\text{Hom}}
\newcommand\cP{\Cal P}

\newcommand\cH{\Cal H}

\newcommand\pa{\partial}

\newcommand\pis{\pi i \sigma}
\newcommand\sd{\,\,{\vartriangleright}\kern -1.0ex{<}\,}
\newcommand\wt{\widetilde}
\newcommand\fg{\frak g}
\newcommand\fk{\frak k}
\newcommand\fp{\frak p}
\newcommand\fs{\frak s}
\newcommand\fh{\frak h}
\newcommand\Cal{\mathcal}

\newcommand\fn{{\frak n}}
\newcommand\fa{{\frak a}}
\newcommand\fm{{\frak m}}
\newcommand\fq{{\frak q}}
\newcommand\CP{{\mathcal P}_n}
\newcommand\Hnm{{\mathbb H}_n \times {\mathbb C}^{(m,n)}}
\newcommand\BD{\mathbb D}
\newcommand\BH{\mathbb H}
\newcommand\CCF{{\mathcal F}_n}
\newcommand\CM{{\mathcal M}}
\newcommand\Gnm{\Gamma_{n,m}}
\newcommand\Cmn{{\mathbb C}^{(m,n)}}
\newcommand\Yd{{{\partial}\over {\partial Y}}}
\newcommand\Vd{{{\partial}\over {\partial V}}}

\newcommand\Ys{Y^{\ast}}
\newcommand\Vs{V^{\ast}}
\newcommand\LO{L_{\Omega}}
\newcommand\fac{{\frak a}_{\mathbb C}^{\ast}}

\newcommand\POB{ {{\partial}\over {\partial{\overline \Omega}}} }
\newcommand\PZB{ {{\partial}\over {\partial{\overline Z}}} }
\newcommand\PX{ {{\partial}\over{\partial X}} }
\newcommand\PY{ {{\partial}\over {\partial Y}} }
\newcommand\PU{ {{\partial}\over{\partial U}} }
\newcommand\PV{ {{\partial}\over{\partial V}} }
\newcommand\PO{ {{\partial}\over{\partial \Omega}} }
\newcommand\PZ{ {{\partial}\over{\partial Z}} }

\begin{section}{{\bf Introduction}}
\setcounter{equation}{0}
\vskip 3mm
Let
$$
G=Sp(4,\BR)=\left\{ M\in GL(4,\BR)\,|\ {}^tMJ_4 M=J_4\,\right\}
$$
be the symplectic group of degree 2, where
$$
J_4=\begin{pmatrix}
  \ 0 & I_2 \\
  -I_2 & 0
\end{pmatrix}.
$$
Let
$$
\BH_2:=\left\{ \Om\in \BC^{(2,2)}\,|\ \Om={}^t\Om,\ {\rm Im}\,\Om>0\,\right\}
$$
be the Siegel upper half plane of degree 2. The $G$ acts on $\BH_2$ transitively by
\begin{equation}\label{(1.1)}
  M\langle \Om\rangle =(A\Om+B)(C\Om+D)^{-1},
\end{equation}
where
$M=\begin{pmatrix}
     A & B \\
    C & D
   \end{pmatrix}\in G$ and $\Om\in\BH_2.$ The stabilizer $K$ of the action (1.1) at
   $i\,I_2$ is given by
$$
K=\left\{ \begin{pmatrix}
     A & B \\
    C & D
   \end{pmatrix}\in G\,\Big|\ A\,{}^t\!A+B\,{}^t\!B=I_2,\ A\,{}^t\!B=B\,{}^t\!A\ \right\}.
$$
$K$ is a maximal compact subgroup of $G$. The mapping
$$
G/K \lrt \BH_2,\qquad gK\longmapsto g\langle iI_2\rangle
$$
is a biholomorphic mapping from the homogeneous space $G/K$ to an Einstein-K{\"a}hler
Hermitian symmetric space $\BH_2$ of dimension three.
\vskip 2mm
We consider two $4\times 4$ real matrices $P$ and $Q$ given by
\begin{equation*}
  P=\begin{pmatrix}
     p & 0 \\
     0 & p
   \end{pmatrix} \quad {\rm and}\quad
   Q=\begin{pmatrix}
     q & 0 \\
     0 & q
   \end{pmatrix},
\end{equation*}
where
\begin{equation*}
  p=\frac{1}{\sqrt{2}}\begin{pmatrix}
     1 & -1 \\
     1 & \ 1
   \end{pmatrix} \quad {\rm and}\quad
   q=\begin{pmatrix}
     0 & 1 \\
     1 & 0
   \end{pmatrix}.
\end{equation*}
Since $q^2=I_2$ and $Q^2=I_4,$ we see that $Q\in G$.
\vskip 3mm
Wang, Wang and Lu (cf.\,\cite{WWL}) introduced the new kind Siegel upper half space
\begin{equation}\label{(1.2)}
  \widehat{\BH}_2:=\left\{ \Om\in\BH_2\,|\ Q\langle \Om \rangle=\Om \,\right\}.
\end{equation}
Then it is easily seen that
\begin{equation}\label{(1.3)}
  \widehat{\BH}_2=\left\{ \Om=
  \begin{pmatrix}
     \tau & z\\
     z & \tau
   \end{pmatrix}
  \in\BH_2\,\big|\ \tau,z\in \BC,\ {\rm Im}\,\tau > |{\rm Im}\,z|\geq 0 \,\right\}.
\end{equation}
$\widehat{\BH}_2$ is a two dimensional complex submanifold of $\BH_2.$

\vskip 3mm
We fix an element
$\Om=\begin{pmatrix}
     \tau & z\\
     z & \tau
    \end{pmatrix}\in \widehat{\BH}_2.$ Let
$$ L_\Om:=\BZ^{(1,2)}\Om+\BZ^{(1,2)}$$
be the lattice in $\BC^2$ determined by $\Om$. We observe that
$$ A_\Om=\BC^2/L_\Om$$
is an abelian surface which is very special.

\vskip 3mm
The aim of this article is to study holomorphic line bundles over an abelian
surface $A_\Omega$ explicitly. This paper is organized as follows. In section 2, we
investigate the new kind Siegel upper space $\widehat{\BH}_2$ in some details.
We also survey some geometric results on $\widehat{\BH}_2$ obtained by the authors
in \cite{WWL}.
In section 3, we study holomorphic line bundles over
the special abelian surface $A_\Om$ explicitly.

\vskip 0.51cm \noindent
{\bf Notations:} \ \ We denote by
$\BR$ and $\BC$ the field of real numbers and the field of complex numbers
respectively. We denote by $\BZ$ and $\BZ^+$ the ring of integers and the set of
all positive integers respectively. The symbol ``:='' means that
the expression on the right is the definition of that on the left.
For two positive integers $k$ and $l$, $F^{(k,l)}$ denotes the set
of all $k\times l$ matrices with entries in a commutative ring $F$.
For any $M\in F^{(k,l)},\
^t\!M$ denotes the transpose of a matrix $M$. $I_n$ denotes the
identity matrix of degree $n$. For a complex matrix $B$,
${\overline B}$ denotes the complex {\it conjugate} of $B$.
For a square matrix $C\in F^{(k,k)}$ of degree $k$,
${\rm Tr}\,(C)$ denotes the trace of $C$.
We denote $T=\{ t\in\BC\,|\,|t|=1\,\}$. $E:={\rm Im}\,H$ denotes the
imaginary part of a Riemann form $H$. $\widehat{A}$ denotes the dual abelian
surface of an abelian surface $A$.  $\mathfrak{P}$ denotes the Poincar{\'e}
bundle over $A\times \widehat{A}$.

\end{section}

\vskip 10mm

\begin{section}{{\bf The new kind Siegel upper half space }}
\setcounter{equation}{0}
\vskip 2mm
Let
\begin{equation}\label{(2.1)}
  \widehat{G}:=\left\{ M\in G\,|\ M\langle \Om\rangle\in \widehat{\BH}_2\ \ {\rm for\ all}
  \ \Om\in \widehat{\BH}_2\,\right\}
\end{equation}
be the subgroup of $G$.
\begin{theorem}\label{thm:2.1}
Let $\varepsilon=\pm 1.$ Then
$$ \widehat{G}=\left\{ M\in Sp(4,\BR)\,|\ MQ=\varepsilon\, QM\,\right\}.$$
Here
$$  Q=\begin{pmatrix}
     q & 0 \\
     0 & q
   \end{pmatrix}
=Q^{-1}\ {\rm with}\
   q=\begin{pmatrix}
     0 & 1 \\
     1 & 0
   \end{pmatrix}=q^{-1}.$$
\end{theorem}
\begin{proof}
See Theorem 1 in \cite[pp.\,4--5]{WWL}.
\end{proof}

\vskip 2mm
Let
\begin{equation*}
  \BD_2:=\left\{ W\in \BC^{(2,2)}\,|\ W=\,{}^t\!W,\ \, I_2-W{\overline W}> 0\,\right\}
\end{equation*}
be the bounded symmetric domain in $\BC^3.$ The Cayley transform $\Phi:\BD_2\lrt \BH_2$
defined  by
\begin{equation}\label{(2.2)}
  \Phi(W):=i(I_2+W)(I_2-W)^{-1},\quad W\in \BD_2
\end{equation}
is a biholomorphic mapping of $\BD_2$ onto $\BH_2$\,(cf.\,\cite[p.\,782]{Y2}).
The inverse mapping $\Psi:\BH_2\lrt\BD_2$ of $\Phi$ is given by
\begin{equation}\label{(2.3)}
  \Psi(\Om):=(\Om-iI_2)(\Om+iI_2)^{-1},\quad \Om\in \BH_2.
\end{equation}
Let
\begin{equation*}
  T=\frac{1}{\sqrt{2}}
  \begin{pmatrix}
     \,I_2 & \ I_2 \\
     iI_2 & -iI_2
   \end{pmatrix}
\end{equation*}
be the $4\times 4$ matrix represented by the Cayley transform $\Phi$. Then
\begin{equation}\label{(2.4)}
  T^{-1}GT=\left\{ \begin{pmatrix}
     \alpha & \beta \\
     \overline{\beta} & \overline{\alpha}\end{pmatrix}
     \in \BC^{(4,4)}\,\big|\
     {}^t\alpha \overline{\alpha}-{}^t\overline{\beta}\beta=I_2,\ \,
     {}^t\alpha \overline{\beta}=\,{}^t\overline{\beta}\alpha\,\right\}.
\end{equation}
Indeed, if $M=\begin{pmatrix}
     A & B \\
     C & D
   \end{pmatrix}\in G$, then
\begin{equation}\label{(2.5)}
  T^{-1}MT=\begin{pmatrix}
     \alpha & \beta \\
     \overline{\beta} & \overline{\alpha}\end{pmatrix},
\end{equation}
where
$$
\alpha=\frac{1}{2}\left\{ (A+D)+i(B-C)\right\}\quad {\rm and}\quad
\beta=\frac{1}{2}\left\{ (A-D)-i(B+C)\right\}.
$$
For brevity, we put $G_*:=T^{-1}GT.$ Then $G_*$ is a ${\it proper}$ subgroup of
$SU(2,2)$, where
$$
SU(2,2):=\{\,h\in\BC^{(4,4)}\,|\ {}^t\!h I_{2,2}\overline{h}=I_{2,2}\,\},\quad
\ I_{2,2}:=\begin{pmatrix}
     I_2 & \ 0 \\
     0 & -I_2\end{pmatrix}.
$$
In fact, since ${}^tTJ_4T=-iJ_4,$ we get
\begin{equation}\label{(2.6)}
  G_*=\{ h\in SU(2,2)],|\ {}^t\!hJ_4h=J_4\,\}=SU(2,2)\cap Sp(4,\BC),
\end{equation}
where
$$Sp(4,\BC)=\{ \alpha\in \BC^{(4,4)}\,|\ {}^t\!\alpha J_4 \alpha=J_4\,\}.$$

\vskip 2mm
Let
\begin{equation*}
  P^+=\left\{ \begin{pmatrix}
                I_2 & Z \\
                0 & I_2
              \end{pmatrix}\,\big|\ Z=\,{}^tZ\in \BC^{(2,2)}\,\right\}
\end{equation*}
be the $P^+$-part of the complexification of $G_*\subset SU(2,2).$
We note that the Harish-Chandra decomposition of an element
$\begin{pmatrix}
     \alpha & \beta \\
     \overline{\beta} & \overline{\alpha}\end{pmatrix}\in G_*$ is
\begin{equation*}
  \begin{pmatrix}
     \alpha & \beta \\
     \overline{\beta} & \overline{\alpha}\end{pmatrix}=
  \begin{pmatrix}
     I_2 & \beta \overline{\alpha}^{-1} \\
     0 & I_2\end{pmatrix}
  \begin{pmatrix}
     \alpha-\beta\overline{\alpha}^{-1}\overline{\beta} & 0 \\
     0 & \overline{\alpha}\end{pmatrix}
 \begin{pmatrix}
     I_2 & 0 \\
     \overline{\alpha}^{-1}\overline{\beta} & I_2\end{pmatrix}.
\end{equation*}
For more detail, we refer to \cite[p.\,155]{K}. Thus the $P^+$-component of
the following element
\begin{equation*}
  \begin{pmatrix}
     \alpha & \beta \\
     \overline{\beta} & \overline{\alpha}\end{pmatrix}\cdot
  \begin{pmatrix}
     I_2 & W \\
     0 & I_2\end{pmatrix},\quad W\in \BD_2
\end{equation*}
of the complexification of $G_*$ is given by
\begin{equation}\label{(2.7)}
  \begin{pmatrix}
     I_2 & (\alpha W+\beta)(\overline{\beta} W+\overline{\alpha})^{-1} \\
     0 & I_2\end{pmatrix}.
\end{equation}
We observe that $\beta \overline{\alpha}^{-1}\in \BD_2.$ We get the
Harish-Chandra embedding of $\BD_2$ into $P^+$ (cf. \cite[p.\,155]{K} or
\cite[pp.\,58--59]{S}). Therefore we see that $G_*$ acts on $\BD_2$ transitively
by
\begin{equation}\label{(2.8)}
  \begin{pmatrix}
     \alpha & \beta \\
     \overline{\beta} & \overline{\alpha}\end{pmatrix}
     \langle W\rangle=(\alpha W+\beta)(\overline{\beta} W+\overline{\alpha})^{-1},
     \quad \begin{pmatrix}
     \alpha & \beta \\
     \overline{\beta} & \overline{\alpha}\end{pmatrix}\in G_*,\ W\in \BD_2.
\end{equation}
The stabilizer of the action (2.8) at the origin $o$ is given by
$$
K_*=\left\{
\begin{pmatrix}
  \alpha & 0 \\
  0 & \overline{\alpha}\end{pmatrix}\in G_* \,\big|\ \alpha\in U(2)\,\right\}.
$$
Thus the homogeneous manifold $G_*/ K_*$ is biholomorphic to $\BD_2$. It is known that
the action \eqref{(1.1)} of $G$ on $\BH_2$ is compatible with
the action \eqref{(2.8)} of $G_*$ on
$\BD_2$ via the Cayley transform $\Phi$. In other words, if $M\in G$ and $W\in \BD_2$, then
\begin{equation}\label{(2.9)}
  M\langle \Phi(W)\rangle=\Phi (M_*\langle W\rangle),
\end{equation}
where $M_*=T^{-1}MT\in G_*.$ The proof can be found in \cite[pp.\,788--789]{Y2}.

\begin{definition}\label{def:2.1}
$$
\widehat{G}_*:=T^{-1}\widehat{G}T\quad {\rm and}\quad \widehat{\BD}_2:=
\Phi^{-1}(\widehat{\BH}_2).
$$
\end{definition}

\begin{lemma}\label{lem:2.1}
Let $\varepsilon=\pm 1.$ Then
$$
 \widehat{G}_*=\left\{ h\in \BC^{(4,4)}\,|\ J_4 [h]=\,{}^th J_4 h=J_4,\ \,
Qh=\varepsilon\,hQ\,\right\}.
$$
\end{lemma}

\begin{proof}
See Formula (16) in \cite[p.\,9]{WWL}.
\end{proof}
According to Lemma \ref{lem:2.1}, we see that
$h=\begin{pmatrix}
     \alpha & \beta \\
     \overline{\beta} & \overline{\alpha}\end{pmatrix}$ with $\alpha,\beta\in\BC^{(2,2)}$
is an element of $\widehat{G}_*$ if and only if
$\alpha\,{}^t\overline{\alpha}-\beta\,{}^t\overline{\beta}=I_2,\
\alpha\,{}^t\overline{\beta}=\overline{\beta}\,{}^t\alpha,\ Qh=\varepsilon\,hQ.$

\begin{lemma}\label{lem:2.2}
$$
\widehat{\BD}_2=\{ W\in \BD_2\,|\ qW=Wq,\ \,I_2-W\overline{W}>0\,\}.
$$
\end{lemma}
\begin{proof}
See Formula (10) in \cite[p.\,7]{WWL}.
\end{proof}

\begin{theorem}\label{thm:2.2}
$\widehat{G}$ acts on $\widehat{\BH}_2$ transitively and
$\widehat{G}_*$ acts on $\widehat{\BD}_2$ transitively.
\end{theorem}
\begin{proof}
See Theorem 2 in \cite[p.\,10]{WWL}.
\end{proof}

\begin{theorem}\label{thm:2.3}
Let ${\sf Bihol}(\widehat{\BH}_2)$ be the group of all biholomorphic mappings of
$\widehat{\BH}_2$ onto $\widehat{\BH}_2$ and ${\sf Bihol}(\widehat{\BD}_2)$ be the group of all biholomorphic mappings of $\widehat{\BD}_2$ onto $\widehat{\BD}_2$. Then
$$
{\sf Bihol}(\widehat{\BH}_2)\cong \widehat{G}/\{\pm I_2,\pm Q\}\quad {\rm and}\quad
{\sf Bihol}(\widehat{\BD}_2)\cong \widehat{G}_*/\{\pm I_2,\pm Q\}.
$$
\end{theorem}
\begin{proof}
See Theorem 3 in \cite[p.\,12]{WWL}.
\end{proof}

\begin{theorem}\label{thm:2.4}
The action of $\widehat{G}$ on $\widehat{\BH}_2$ is compatible with the
action of $\widehat{G}_*$ on $\widehat{\BD}_2$ via the Cayley transform $\Phi$.
More precisely, if $M\in \widehat{G}$ and $W\in \widehat{\BD}_2$, then
\begin{equation}\label{(2.10)}
  M\langle \Phi(W)\rangle=\Phi (M_*\langle W\rangle),
\end{equation}
where $M_*=T^{-1}MT\in \widehat{G}_*.$
\end{theorem}
\begin{proof}
It follows immediately from Theorem \ref{thm:2.2} and Formula \eqref{(2.9)}.
\end{proof}

\vskip 2mm
According to Siegel's work \cite{Si} or \cite[Theorem 6, pp.\,20--21]{WWL},
the symplectic metric on $\widehat{\BH}_2$ is given by
$$ ds^2={\rm Tr}(Y^{-1}d\Om Y^{-1} d\overline{\Om}),$$
where $\Om=X+iY\in \widehat{\BH}_2$ with $X,Y\in\BR^{(2,2)}$ and $Y>0.$
Clearly $ds^2$ is invariant under the action of $\widehat{G}$ on $\widehat{\BH}_2$.
If we take a coordinate $\Om=\begin{pmatrix}
                               \tau & z \\
                               z & \tau
                             \end{pmatrix}\in \widehat{\BH}_2$ with
$\tau=x+iy,\ z=u+iv,\ x,y,u,v\in\BR,$
then $Y^{-1}=(y^2-v^2)^{-1}\begin{pmatrix}
                                \ y & -v \\
                                -v & \ y
                              \end{pmatrix}$ and
\begin{equation*}
  Y^{-1}d\Om=\frac{1}{y^2-v^2}\begin{pmatrix}
                                y\,d\tau-v\,dz & y\,dz-v\,d\tau \\
                                y\,dz-v\,d\tau & y\,d\tau-v\,dz
                              \end{pmatrix}.
\end{equation*}
By a direct computation, we get
\begin{eqnarray*}
  ds^2\!\!\! &=&\!\!\!{\rm Tr}(Y^{-1}d\Om\, Y^{-1}d\overline{\Om}) \\
  {}\!\!\! &=& \!\!\!\frac{2\,(y^2+v^2)}{(y^2-v^2)^2}(d\tau\,d\overline{\tau}+dz\,\overline{z})
  -\frac{4\,yv}{(y^2-v^2)^2}(d\tau\,d\overline{z} +d\overline{\tau}\,dz)\\
  {}\!\!\! &=&\!\!\! \frac{2\,(y^2+v^2)}{(y^2-v^2)^2}
  (dx^2+dy^2+du^2+dv^2)-\frac{8\,yv}{(y^2-v^2)^2}(dx\,du+dy\,dv).
\end{eqnarray*}
We put
\begin{equation*}
  \frac{\partial}{\partial\Om}=
  \begin{pmatrix}
  \,\frac{\partial}{\partial\tau} & \frac{1}{2}\frac{\partial}{\partial z} \\
  \frac{1}{2}\frac{\partial}{\partial z} & \,\frac{\partial}{\partial\tau}
  \end{pmatrix}
  \quad {\rm and}\quad
   \frac{\partial}{\partial\overline{\Om}}=
  \begin{pmatrix}
  \,\frac{\partial}{\partial\overline{\tau}} & \frac{1}{2}\frac{\partial}{\partial\overline{z}} \\
  \frac{1}{2}\frac{\partial}{\partial\overline{z}} & \,\frac{\partial}{\partial\overline{\tau}}
  \end{pmatrix}.
\end{equation*}
H. Maass \cite{M} proved that the differential operator
\begin{equation*}
\Delta=\,{4\over A}\,{\rm Tr} \left(Y\,
{}^{{}^{{}^{{}^\text{\scriptsize $t$}}}}\!\!\!
\left(Y\POB\right)\PO\right)
\end{equation*}
is the Laplace operator $\Delta$ of $ds^2$.
Since
\begin{equation*}
  Y\frac{\partial}{\partial\overline{\Om}}
=\begin{pmatrix}
y\,\frac{\partial}{\partial\overline{\tau}} + \frac{1}{2}\,v\,\frac{\partial}{\partial\overline{z}}
& v\,\frac{\partial}{\partial\overline{\tau}} + \frac{1}{2}\,y\,\frac{\partial}{\partial\overline{z}}\\
  v\,\frac{\partial}{\partial\overline{\tau}} + \frac{1}{2}\,y\,\frac{\partial}{\partial\overline{z}}
& y\,\frac{\partial}{\partial\overline{\tau}} + \frac{1}{2}\,v\,\frac{\partial}{\partial\overline{z}}
\end{pmatrix}
={}^{{}^{{}^{{}^\text{\scriptsize $t$}}}}\!\!\!
\left(Y\POB\right),
\end{equation*}
we get
\begin{eqnarray*}
  \Delta\!\!\! &=&\!\!\! 4\,(y^2+v^2)\,\frac{\partial^2}{\partial{\tau}\partial\overline{\tau}}
  +(y^2+v^2)\,\frac{\partial^2}{\partial{z}\partial\overline{z}}
  +4\,yv\,\left( \frac{\partial^2}{\partial{\tau}\partial\overline{z}} +
  \frac{\partial^2}{\partial{z}\partial\overline{\tau}} \right) \\
  {}\!\!\! &=&\!\!\! (y^2+v^2)
  \,\left( \frac{\partial^2}{\partial x^2}+\frac{\partial^2}{\partial y^2}\right)
  +\frac{1}{4}\,(y^2+v^2)\,\left( \frac{\partial^2}{\partial u^2}+\frac{\partial^2}{\partial v^2}\right)
  +yv\left( \frac{\partial^2}{\partial x \partial u} +\frac{\partial^2}{\partial y \partial v}\right).
\end{eqnarray*}
And
\begin{equation}\label{(2.11)}
  dv(\Om)=\frac{4}{(y+v)^2(y-v)^2}\,dxdydudv
\end{equation}
is a volume element invariant under the action \eqref{(1.1)} of $\widehat{G}$ on
$\widehat{\BH}_2$ (cf. \cite[Theorem 10, pp.\,40--43]{WWL}).

\begin{theorem}\label{thm:2.5}
Let $\Om_1=\begin{pmatrix}
             \tau_1 & z_1 \\
             z_1 & \tau_1
           \end{pmatrix}$ and
$\Om_2=\begin{pmatrix}
             \tau_2 & z_2 \\
             z_2 & \tau_2
           \end{pmatrix}$ be any two points in $\widehat{\BH}_2$. We write
\begin{eqnarray*}
  \tau_1+z_1\!\!&=&\!\! x_1+iy_1,\ \,\tau_1-z_1= x_2+iy_2,\ \,x_1,y_1,x_2,y_2\in \BR,\\
  \tau_2+z_2\!\!&=&\!\! u_1+iv_1,\ \,\tau_2-z_2= u_2+iv_2,\ \,u_1,v_1,u_2,v_2\in \BR,\\
  A \!\!&=&\!\! \frac{y_1^2+v_1^2+(x_1-u_1)^2}{y_1v_1}\geq 2,\ \,
  B= \frac{y_2^2+v_2^2+(x_2-u_2)^2}{y_2v_2}\geq 2,  \\
  \lambda\!\! &=&\!\! \frac{A+\sqrt{A^2-4}}{2}\geq 1 \quad {\rm and}\quad
  \mu= \frac{B+\sqrt{B^2-4}}{2}\geq 1.
\end{eqnarray*}
Then the distance $\rho(\Om_1,\Om_2)$ between $\Om_1$ and $\Om_2$ is given by
$$
\rho(\Om_1,\Om_2)=\left\{ (\log \lambda)^2+(\log \mu)^2\right\}^{1/2}.
$$
\end{theorem}
\begin{proof}
See Theorem 8 in \cite[pp.\,33--34]{WWL}.
\end{proof}

\begin{theorem}\label{thm:2.6}
Let $\Om_1=\begin{pmatrix}
             \tau_1 & z_1 \\
             z_1 & \tau_1
           \end{pmatrix}$ and
$\Om_2=\begin{pmatrix}
             \tau_2 & z_2 \\
             z_2 & \tau_2
           \end{pmatrix}$ be any two points in $\widehat{\BH}_2$. We write
\begin{eqnarray*}
  \tau_1+z_1\!\!&=&\!\! x_1+iy_1,\ \,\tau_1-z_1= x_2+iy_2,\ \,x_1,y_1,x_2,y_2\in \BR,\\
  \tau_2+z_2\!\!&=&\!\! u_1+iv_1,\ \,\tau_2-z_2= u_2+iv_2,\ \,u_1,v_1,u_2,v_2\in \BR,\\
  A \!\!&=&\!\! \frac{y_1^2+v_1^2+(x_1-u_1)^2}{y_1v_1}\geq 2,\ \,
  B= \frac{y_2^2+v_2^2+(x_2-u_2)^2}{y_2v_2}\geq 2,  \\
  \lambda\!\! &=&\!\! \frac{A+\sqrt{A^2-4}}{2}\geq 1 \quad {\rm and}\quad
  \mu= \frac{B+\sqrt{B^2-4}}{2}\geq 1.
\end{eqnarray*}
Then the geodesic $\gamma(s)$ joining $\Om_1$ and $\Om_2$ can be expressed as
\begin{equation}\label{(2.11)}
 \gamma(s)=\begin{pmatrix}
             \tau(s) & z(s) \\
             z(s) & \tau(s)
           \end{pmatrix},
\end{equation}
\begin{eqnarray*}
  \tau(s): &=& \frac{x_1+x_2}{2}+\frac{y_1}{2} R(s)+\frac{y_2}{2} T(s),\\
  z(s): &=& \frac{x_1-x_2}{2}+\frac{y_1}{2} R(s)-\frac{y_2}{2} T(s),
\end{eqnarray*}
where
\begin{eqnarray*}
  R(s) &=& \frac{\lambda (u_1-x_1)\left( \lambda^{2s/s_0}-1\right)
  c(\lambda)+i\,v_1 \lambda^{s/s_0}} {(\lambda y_1-v_1)\left( \lambda^{2s/s_0}-1\right)
  c(\lambda)+v_1}, \\
  T(s) &=& \frac{\mu (u_2-x_2)\left( \mu^{2s/s_0}-1\right)
  c(\mu)+i\,v_2 \mu^{s/s_0}}
  {(\mu y_2-v_2)\left( \mu^{2s/s_0}-1\right) c(\mu)+v_2}.
\end{eqnarray*}
\vskip 2mm\noindent
And $c(x)$ is the function defined on the interval $[1,\infty]$ by
\begin{equation*}
  c(x)=\begin{cases}
         \frac{1}{x^2-1} & \mbox{if } x>1 \\
         \ \ 1 & \mbox{if}\ x=1
       \end{cases}
\end{equation*}
and the parameter $s$ (with $0\leq s < s_0)$ denotes the arc length of the arc initiated at
$\Om_1$ in the geodesic.
\end{theorem}
\begin{proof}
See Theorem 9 in \cite[pp.\,40--43]{WWL}.
\end{proof}

\vskip 3mm
Now we investigate the group $\widehat{G}$.
\begin{proposition}\label{prop:2.1}
Let $\varepsilon=\pm 1.$ Then we have the following\,:\\
\noindent
${\rm (1)}$
\begin{equation*}
  \widehat{G}=\left\{ \begin{pmatrix}
              a_1 & a_2 & b_1 & b_2 \\
              \varepsilon a_2 & \varepsilon a_1 & \varepsilon b_2 & \varepsilon b_1 \\
              c_1 & c_2 & d_1 & d_2 \\
              \varepsilon c_2 & \varepsilon c_1 & \varepsilon d_2 & \varepsilon d_1
            \end{pmatrix}\in \BR^{(4,4)} \Bigg|\
            \begin{pmatrix}
             a_1d_1+a_2d_2-b_1c_1-b_2c_2= 1,  \\
              a_1d_2+a_2d_1-b_2c_1-b_1c_2=0
            \end{pmatrix}\ \right\}.
\end{equation*}
\noindent
${\rm (2)}$ We put
\begin{equation*}
  \widehat{G}_+=\left\{ \begin{pmatrix}
              a_1 & a_2 & b_1 & b_2 \\
               a_2 &  a_1 &  b_2 &  b_1 \\
              c_1 & c_2 & d_1 & d_2 \\
               c_2 &  c_1 &  d_2 &  d_1
            \end{pmatrix} \in\BR^{(4,4)}\Bigg|\  a_1d_1+a_2d_2-b_1c_1-b_2c_2=1,\
            a_1d_2+a_2d_1-b_2c_1-b_1c_2=0\ \right\}.
\end{equation*}
Then $\widehat{G}_+$ is a proper subgroup of $\widehat{G}$ and
$\dim \widehat{G}=\dim \widehat{G}_+=6.$ The Lie algebra $\widehat{\mathfrak g}$
of $\widehat{G}$ is given by
$$
\widehat{\mathfrak g}
=\left\{
\begin{pmatrix}
 X & \, Y \\
 Z & -{}^tX
\end{pmatrix}\in \BR^{(4,4)} \big|\ X,Y,Z\in \mathbb A\,\right\},
$$
where
$$
\mathbb A:=\left\{
\begin{pmatrix}
  x &  y \\
  y &  x
\end{pmatrix}\in \BR^{(2,2)}\,\,\right\}.
$$
The Lie algebra $\widehat{\mathfrak g}_+$ of $\widehat{G}_+$ is equal to
$\widehat{\mathfrak g}$.\\
\vskip 2mm
\noindent
${\rm (3)}$ The stabilizer $\widehat{K}$ of the action of $\widehat{G}$ at $iI_2$ is given by
\begin{equation*}
  \widehat{K}=\left\{ \begin{pmatrix}
              \ a_1 & \ a_2 & b_1 & b_2 \\
               \varepsilon a_2 & \varepsilon a_1 & \varepsilon b_2 & \varepsilon b_1 \\
              -b_1 & -b_2 & a_1 & a_2 \\
               -\varepsilon b_2 & -\varepsilon b_1 & \varepsilon a_2 & \varepsilon a_1
            \end{pmatrix} \in\BR^{(4,4)}\Bigg|\  a_1^2+a_2^2+b_1^2+b_2^2=1,\
            a_1a_2+b_1b_2=0\ \right\}.
\end{equation*}
The stabilizer $\widehat{K}_+$ of the action of $\widehat{G}_+$ at $iI_2$ is given by
\begin{equation*}
  \widehat{K}_+=\left\{ \begin{pmatrix}
              \ a_1 & \ a_2 & b_1 & b_2 \\
               \ a_2 & \ a_1 &  b_2 &  b_1 \\
              -b_1 & -b_2 & a_1 & a_2 \\
               -b_2 &  -b_1 &  a_2 &  a_1
            \end{pmatrix} \in\BR^{(4,4)}\Bigg|\  a_1^2+a_2^2+b_1^2+b_2^2=1,\
            a_1a_2+b_1b_2=0\ \right\}.
\end{equation*}
Thus $\dim \widehat{K}=\dim \widehat{K}_+=2$. Let $\widehat{\mathfrak k}$
\ (resp. $\widehat{\mathfrak k}_+)$ be the Lie algebra of $\widehat{K}$ (resp. $\widehat{K}_+$).
Then
\begin{equation*}
  \widehat{\mathfrak k}=\widehat{\mathfrak k}_+=
  \left\{ \begin{pmatrix}
            \ 0 & Z \\
            -Z & 0
          \end{pmatrix}
  \in \BR^{(4,4)}\,\Big|\ Z=\begin{pmatrix}
            x & y \\
            y & x
          \end{pmatrix}\in \BR^{(2,2)}\,\right\}.
\end{equation*}
\end{proposition}
\begin{proof}
The proof of (1) follows from Theorem \ref{thm:2.1}. The proof of (2) is given as follows.
We know that the Lie algebra $\mathfrak{sp}(4,\BR)$ of $Sp(4,\BR)$ is
\begin{equation*}
  \mathfrak{sp}(4,\BR)=\left\{ \begin{pmatrix}
                                 X & \ Y \\
                                 Z & -{}^t\!X
                               \end{pmatrix}\in \BR^{(4,4)}\,\big|\
                               X,Y,Z\in\BR^{(2,2)},\ Y={}^tY,\ Z={}^tZ\ \right\}.
\end{equation*}
From (1) and the relation $MQ=\varepsilon\,QM$ for all $M\in \widehat{G}$, we see that
$$
\widehat{\mathfrak g}
=\left\{
\begin{pmatrix}
 X & \, Y \\
 Z & -{}^tX
\end{pmatrix}\in \BR^{(4,4)} \big|\ X,Y,Z\in \mathbb A\,\right\},
$$
where
$$
\mathbb A:=\left\{
\begin{pmatrix}
  x &  y \\
  y &  x
\end{pmatrix}\in \BR^{(2,2)}\,\,\right\}.
$$
Since $\dim \widehat{G}=\dim \widehat{G}_+=6,$ we have
$\widehat{\mathfrak g}=\widehat{\mathfrak g}_+.$
The proof of (3) follows from (a), (b) and a direct calculation.
\end{proof}

\begin{lemma}\label{lem:2.3}
We put
\begin{equation*}
  \widehat{G}_-=\left\{ \begin{pmatrix}
              \,a_1 & \,a_2 & \,b_1 & \,b_2 \\
               -a_2 &  -a_1 &  -b_2 &  -b_1 \\
              \,c_1 & \,c_2 & \,d_1 & \,d_2 \\
               -c_2 &  -c_1 &  -d_2 &  -d_1
            \end{pmatrix}\in \widehat{G} \, \right\}.
\end{equation*}
Then $\widehat{G}_-$ is not a subgroup of $\widehat{G}$ and we have the following relations:
$$
\widehat{G}=\widehat{G}_+ \cup \widehat{G}_-,\quad
\widehat{G}_+\!\cdot\! \widehat{G}_+\subset \widehat{G}_+,
\quad \widehat{G}_+\!\cdot\! \widehat{G}_-\subset \widehat{G}_-,\quad
\widehat{G}_-\!\cdot\! \widehat{G}_+\subset \widehat{G}_-,
\quad \widehat{G}_-\!\cdot\! \widehat{G}_-\subset \widehat{G}_+.
$$
\end{lemma}
\begin{proof}
The proof follows from an easy direct computation.
\end{proof}

\begin{lemma}\label{lem:2.4}
Let $M=\begin{pmatrix}
         A & B \\
         C & D
       \end{pmatrix}$ be an element of $\widehat{G}$ with
       $A=\begin{pmatrix}
         a_1 & a_2 \\
         \varepsilon a_2 & \varepsilon a_1
       \end{pmatrix}, \
       B=\begin{pmatrix}
         b_1 & b_2 \\
         \varepsilon b_2 & \varepsilon b_1
       \end{pmatrix},\
       C=\begin{pmatrix}
         c_1 & c_2 \\
         \varepsilon c_2 & \varepsilon c_1
       \end{pmatrix}$  and
       $D=\begin{pmatrix}
         d_1 & d_2 \\
         \varepsilon d_2 & \varepsilon d_1
       \end{pmatrix}\in \BR^{(2,2)}$. Here $\varepsilon=\pm 1.$
       We put
\begin{equation}\label{(2.13)}
  M_1:=\begin{pmatrix}
         a_1+a_2 & b_1+b_2 \\
         c_1+c_2 & d_1+d_2
       \end{pmatrix}\quad {\rm and}\quad
 M_2:=\begin{pmatrix}
         a_1-a_2 & b_1-b_2 \\
         c_1-c_2 & d_1-d_2
       \end{pmatrix}\in \BR^{(2,2)}.
\end{equation}
Then $\det M_1=\det M_2=1$, that is, $M_1,M_2\in SL(2,\BR).$
\end{lemma}
\begin{proof}
According to Proposition \ref{prop:2.1}, if $M=\begin{pmatrix}
         A & B \\
         C & D
       \end{pmatrix}$ is an element of $\widehat{G}$,
\begin{equation}\label{(2.14)}
  a_1d_1+a_2d_2-b_1c_1-b_2c_2=1\quad {\rm and}\quad
            a_1d_2+a_2d_1-b_2c_1-b_1c_2=0.
\end{equation}
Then by Formula (2.14)
\begin{eqnarray*}
  \det M_1\!\! &=&\!\! (a_1+a_2)(d_1+d_2)-(b_1+b_2)(c_1+c_2) \\
  {}\!\! &=&\!\! (a_1d_1+a_2d_2-b_1c_1-b_2c_2)+(a_1d_2+a_2d_1-b_2c_1-b_1c_2)=1
\end{eqnarray*}
and
\begin{eqnarray*}
  \det M_2\!\! &=&\!\! (a_1-a_2)(d_1-d_2)-(b_1-b_2)(c_1-c_2) \\
  {}\!\! &=&\!\! (a_1d_1+a_2d_2-b_1c_1-b_2c_2)-(a_1d_2+a_2d_1-b_2c_1-b_1c_2)=1.
\end{eqnarray*}

\end{proof}

\begin{lemma}\label{lem:2.5}
Let $M, A,B,C,D,M_1 $ and $M_2$ be the matrices given in Lemma \ref{lem:2.4}.
If $M=\begin{pmatrix}
         A & B \\
         C & D
       \end{pmatrix}\in \widehat{G}$ and
       $\Om=\begin{pmatrix}
         \tau & z \\
         z & \tau
       \end{pmatrix}\in \widehat{\BH}_2$, then we have
\begin{equation}\label{(2.15)}
  M\langle \Om \rangle =
       \begin{pmatrix}
         \tau_* & z_* \\
         z_* & \tau_*
       \end{pmatrix},
\end{equation}
where
\begin{equation*}
  \tau_*=\frac{M_1\langle \tau+z\rangle+M_2\langle \tau-z\rangle}{2} \quad
  {\rm and}\quad z_*=
  \frac{\varepsilon\left( M_1\langle \tau+z\rangle-M_2\langle \tau-z\rangle\right)}{2}.
\end{equation*}
\end{lemma}

\begin{proof}
The proof can be found in \cite[Lemma 6, pp.\,29--30]{WWL}.
\end{proof}

\begin{lemma}\label{lem:2.6}
Let ${\mathbb B}:\widehat{\mathfrak g}\times \widehat{\mathfrak g}\lrt \BR$ be the
Killing form of $\widehat{\mathfrak g}$ defined by
\begin{equation*}
  {\mathbb B}(X,Y):={\rm Tr}({\rm ad}(X)\,{\rm ad}(Y)),\quad X,Y\in \widehat{\mathfrak g},
\end{equation*}
where ${\rm ad}:\widehat{\mathfrak g}\lrt {\rm End}(\widehat{\mathfrak g})$ denote
the adjoint representation of $\widehat{\mathfrak g}$. Then $\mathbb{B}$ is
degenerate. Therefore $\widehat{G}$ is not a semisimple Lie group.
\end{lemma}

\begin{proof}
We can show that if $X,Y\in \widehat{\mathfrak g}$ with
\begin{equation*}
  X=\begin{pmatrix}
      a_1 & a_2 & \ b_1 & \ b_2 \\
      a_2 & a_1 & \ b_2 & \ b_1 \\
      c_1 & c_2 & -a_1 & -a_2 \\
      c_2 & c_1 & -a_2 & -a_1
    \end{pmatrix}\quad {\rm and}\quad
  Y=\begin{pmatrix}
      p_1 & p_2 & \ q_1 & \ q_2 \\
      p_2 & p_1 & \ q_2 & \ q_1 \\
      r_1 & r_2 & -p_1 & -p_2 \\
      r_2 & r_1 & -p_2 & -p_1
    \end{pmatrix},
\end{equation*}
then we get
$$ {\mathbb B}(X,Y)= c\,{\rm Tr}(XY)= 2\,c\,(2\,a_1p_1+2\,a_2p_2 +b_1r_1+b_2r_2+c_1 q_1+c_2q_2),$$
where $c$ a nonzero real constant. If we take $X=Y$, then
$${\mathbb B}(X,X)= 4\,c\,(a_1^2+a_2^2+b_1c_1+b_2c_2).$$
Clearly ${\mathbb B}$ is degenerate. Therefore $\widehat{\mathfrak g}$ is not a
semisimple Lie algebra.
\end{proof}

\begin{theorem}\label{thm:2.7}
$\widehat{G}_+$ acts on $\widehat{\BH}_2$ transitively.
\end{theorem}
\begin{proof}
Let $\widehat{G},\ \widehat{G}_+,\ \widehat{K}$ and $\widehat{K}_+$ be the subsets of
$Sp(4,\BR)$ defined in Proposition \ref{prop:2.1}.
According to Theorem \ref{thm:2.2}, we know that $\widehat{G}$ acts on $\widehat{\BH}_2$ transitively.
Thus the homogeneous manifold $\widehat{G}/\widehat{K}$ is biholomorphic to the complex manifold
$\widehat{\BH}_2$ via
\begin{equation*}
  \widehat{G}/\widehat{K}\lrt \widehat{\BH}_2,\qquad g\widehat{K}\longmapsto
  g\langle iI_2\rangle,\ \,g\in \widehat{G}.
\end{equation*}
Here $\widehat{K}$ is the stabilizer of the action of $\widehat{G}$ on $\widehat{\BH}_2$
given by Proposition \ref{prop:2.1} (3). We claim that the following mapping
\begin{equation*}
 \Theta : \widehat{G}_+/\widehat{K}_+\lrt \widehat{G}/\widehat{K},
  \qquad g\widehat{K}_+\longmapsto g\widehat{K},\ \,g\in \widehat{G}_+
\end{equation*}
is bijective. Suppose $g_1\widehat{K}=g_2\widehat{K}$ for $g_1,g_2\in \widehat{G}_+.$
Then $g_1^{-1}g_2\in \widehat{K}$ and so
$g_1^{-1}g_2\in \widehat{K}\cap \widehat{G}_+=\widehat{K}_+.$ Since
$g_1\widehat{K}_+=g_2\widehat{K}_+,\ \Theta$ is injective.
For the surjectivity of $\Theta$, for any $g\widehat{K}\in \widehat{G}/\widehat{K}$
with $g\in \widehat{G}$, we have to find $g_*\in \widehat{G}_+$ such that
$\Theta (g_*\widehat{K}_+)=g\widehat{K}.$ If $g\in \widehat{G}_+,$ we take $g_*=g.$
If $g\in \widehat{G}_-=\widehat{G}\setminus\widehat{G}_+,$ we take $g_*=gk$ with
\begin{equation*}
  k=\begin{pmatrix}
              \ a_1 & \ a_2 & b_1 & b_2 \\
               - a_2 & - a_1 & - b_2 & -b_1 \\
              -b_1 & -b_2 & a_1 & a_2 \\
               \ b_2 & \ b_1 & - a_2 & - a_1
            \end{pmatrix}\in \widehat{K}\setminus \widehat{K}_+\subset \widehat{G}_-.
\end{equation*}
Since $g_*=gk\in \widehat{G}_+$ by Lemma \ref{lem:2.3}, we have
$$\Theta(g_*\widehat{K}_+)=g_*\widehat{K}=gk\widehat{K}=g\widehat{K}.$$
Thus $\Theta$ is surjective. Consequently $\Theta$ is bijective and indeed
biholomorphic. Since $\widehat{G}$ acts on $\widehat{\BH}_2$ transitively by
Theorem \ref{thm:2.2}, $\widehat{G}_+$ acts on $\widehat{\BH}_2$ transitively.
\end{proof}

\end{section}

\vskip 10mm

\begin{section}{{\bf Line bundles over a special abelian surface}}
\setcounter{equation}{0}
\vskip 2mm
Let $A=\BC^2/L$ be an abelian surface, where $L$ is a lattice of rank $4$ in $\BC^2$. Then $A$
is a two-dimensional projective variety which is a connected, compact commutative complex
Lie group. For an integer $k$ with $0\leq k\leq 4$, we denote by $H^k(A)$ the space of all
harmonic $k$-forms on $A$. Since $A\cong (S^1)^4$ topologically,
$\dim_\BC H^k(A)={4 \choose k}$ with $0\leq k\leq 4$. The Hodge decomposition of $A$ is given
\begin{equation*}
  H^2(A,\BC)\cong H^0(A,\Om^2)\oplus H^1(A,\Om^1)\oplus H^2(A,\Om^0),
\end{equation*}
where $\Om^p$ is the sheaf of holomorphic $p$-forms on $A$. We have canonical isomorphisms
\begin{equation*}
  H^r(A,\BZ)\cong \left\{
  {\rm group\ of\ alternating}\ r\!\!-\!{\rm forms}\ L\times\cdots\times L\lrt \BZ\,\right\}.
\end{equation*}
See \cite[pp.\,3--4]{MU} for detail.

\begin{definition}\label{def:3.1}
Let $A=\BC^2/L$ be an abelian surface with a lattice $L$. A Hermitian form $H$ on $\BC^2$
is called a {\sf Riemann\ form} for $A$ if it satisfies the following conditions
${\rm (R1)\ and\ (R2)}$\,:
\vskip 2mm \noindent
${\rm (R1)}$ $H$ is nondegenerate\,;
\vskip 2mm \noindent
${\rm (R2)}$ The imaginary part $E={\rm Im}\,H$ of $H$ is integral valued on $L\times L$.
\end{definition}

\vskip 2mm
Let $A=\BC^2/L$ be an abelian surface and let $H$ be an Riemann form for $A$. A map
$\chi:L\lrt T$ is said to be a ${\sf semi\!\!-\!\!character}$ of $L$ with respect to
$E={\rm Im}\,H$ if it satisfies the following condition
\begin{equation}\label{(3.1)}
  \chi(\alpha+\beta)=\chi(\alpha)\,\chi(\beta) \exp\left\{i\,\pi\,E(\alpha,\beta)\right\},\quad
  \alpha,\beta\in L.
\end{equation}
For any Riemann form $H$ on $\BC^2$, a semi-character $\chi_H$ satisfying the condition
(3.1) always exists uniquely (cf.\,\cite[p.\,20]{MU}).
The mapping $J_{H,\chi_H}:L\times \BC^2\lrt GL(1,\BC)$ defined by
\begin{equation}\label{(3.2)}
  J_{H,\chi_H}(\alpha,z):=\chi_H(\alpha)\exp\left\{\pi H(z,\alpha)+ \frac{\pi}{2}
  H(\alpha,\alpha)\right\},\quad \alpha\in L,\ z\in \BC^2
\end{equation}
satisfies the following equation
\begin{equation*}
  J_{H,\chi_H}(\alpha+\beta,z)= J_{H,\chi_H}(\alpha,\beta+z) J_{H,\chi_H}(\beta,z)
\end{equation*}
for all $\alpha,\beta\in L$ and $z\in \BC^2$. In other words, $J_{H,\chi_H}$ is
an automorphic factor of rank one with respect to the lattice $L$ in $\BC^2$.
Let $L(H,\chi_H)$ be the holomorphic line bundle over $A$ defined by the automorphic factor
$J_{H,\chi_H}$ for $L$. More precisely, $L(H,\chi_H)$ is expressed as follows.
We first observe that the lattice $L$ acts on $\BC^2\times \BC$ by
\begin{equation}\label{(3.3)}
  \rho_\alpha (z,\xi):=(z+\alpha,J_{H,\chi_H}(\alpha,z)\,\xi),
\end{equation}
where $\al\in L,\ z\in\BC^2$ and $\xi\in \BC$.
$L(H,\chi_H)$ is the quotient of the trivial line bundle $\BC^2\times\BC$ for the
action of $L$ given by \eqref{(3.3)}.
The global sections of $L(H,\chi_H)$ are essentially the global sections $\theta$ of the trivial
line bundle $\BC^2\times \BC$ over $\BC^2$ (i.e. holomorphic functions $\theta$ on $\BC^2$)
which are invariant under the action \eqref{(3.3)} of $L$ on $\BC^2\times\BC$, that is,
which satisfies the functional equation
\begin{equation}\label{(3.4)}
  \theta(z+\al)=J_{H,\chi_H}(\alpha,z)\,\theta(z),\quad z\in\BC^2,\ \al\in L.
\end{equation}
In other words, a function $\theta$ is a theta function for $H$ and $\chi_H$.
We can say that the cohomology group $H^0(A,L(H,\chi_H))$ is the space of theta
functions for $H$ and $\chi_H$. According to \cite[Proposition, p.\,26]{MU},
we have
\begin{equation}\label{(3.5)}
  \dim\,H^0(A,L(H,\chi_H))=\sqrt{\det E}.
\end{equation}
According to Theorem of Lefschetz (\cite[pp.\,29--30]{MU}, if $H$ is a positive definite
Riemann form for $A$, the line bundle $L(H,\chi_H)$ is ample, that is, the space of
holomorphic sections $L(H,\chi_H)^{\otimes n}\ (n\geq 3)$ give an embedding of $A$ as a
closed complex submanifold into a projective space $\mathbb{P}^N$. Therefore $A$ is a
smooth projective variety of dimension 2.

\vskip 3mm
Without difficulty we can show that there is an isomorphism of holomorphic line
bundles over $A$
\begin{equation*}
  L(H_1,\chi_1)\otimes  L(H_2,\chi_2)\cong L(H_1+H_2,\chi_1\chi_2).
\end{equation*}
Therefore the set of all $L(H,\chi)$ forms a group under tensor product.
We note that two holomorphic line bundles $L(H,\chi)$ and $L(\widetilde{H},\widetilde{\chi})$
over $A$ are isomorphic if and only if their automorphic factors $J_{H,\chi}$ and
$J_{\widetilde{H},\widetilde{\chi}}$ are holomorphically equivalent, i.e., there exists
a holomorphic mapping $h:\BC^2\lrt GL(1,\BC)$ such that
\begin{equation*}
  h(z+\alpha)J_{H,\chi}(\alpha,z)=h(z)J_{\widetilde{H},\widetilde{\chi}}(\alpha,z)
\end{equation*}
for all $\alpha\in L$ and $z\in\BC^2$.
Let ${\rm Pic}^0(A)$ be the group of holomorphic line bundles over an abelian surface $A$
with zero Chern class. In fact, ${\rm Pic}^0(A)\cong {\rm Hom}(L,T)$ is an abelian surface
called the ${\sf dual\ abelian\ surface}$. From now on, we denote $\widehat{A}:={\rm Pic}^0(A).$

\vskip 3mm
Now we have the famous theorem.
\vskip 2mm\noindent
{\bf Theorem\ of\ Appell}-{\bf Humbert}. {\it Any holomorphic line bundle over an abelian surface
$A$ is isomorphic to the line bundle $L(H,\chi)$ for a uniquely determined
Riemann form $H$ for $A$ and a uniquely determined semi-character $\chi$.}
\begin{proof}
The proof of the above theorem can be found in \cite[pp.\,20--22]{MU}.
\end{proof}

\begin{definition}\label{Definition 3.2}
Let $A=\BC^2/L$ be an abelian surface.
\vskip 2mm \noindent
{\rm (1)} For an element $x\in A$, we define the ${\sf translation}\ T_x$ by $x$ by
\begin{equation*}
  T_x:A\lrt A, \quad y\longmapsto x+y \quad {\rm for\ all}\ y\in A.
\end{equation*}
\vskip 2mm \noindent
{\rm (2)} A holomorphic line bundle $F$ over $A$ is said to be ${\sf homogeneous}$ if
$T_x^*(F)\cong F$ for all \\
\indent
\ \,\,$x\in A.$
\vskip 2mm \noindent
{\rm (3)} Let $F$ be a holomorphic line bundle over $A$. We set
\begin{equation*}
  K(F):=\{ x\in A\,|\ T_x^*(F)\cong F\,\}.
\end{equation*}
\ \ \ \
We define the mapping $\phi_F:A\lrt \widehat{A}$ by
\begin{equation}\label{(3.6)}
  \phi_F(x):=T_x^*(F)\otimes F^{-1},\quad x\in A.
\end{equation}
\end{definition}

\begin{remark}\label{Remark 3.1} (1) $K(F)$ is a Zariski-closed subgroup of
$A$ (${\rm cf.}$\,see \cite[Proposition, p.\,60]{MU}).
\vskip 2mm \noindent
{\rm (2)} If $F$ is an ample holomorphic line bundle over $A$, then $\phi_F$
is an isogeny, i.e., $\phi_F$ is a surjective homomorphism with the finite
kernel $K(F)$.
\end{remark}

\vskip 3mm\noindent
{\bf Theorem of Square.} {\it Let $F$ be an ample holomorphic line bundle over an
abelian surface $A$. Then for any $x,y\in A$, }
\begin{equation*}
  T^*_{x+y}(F)\otimes F \cong T_x^*((F)\otimes T_y^*(F).
\end{equation*}
\begin{proof}
See \cite[Corollary 4, p.\,59]{MU}.
\end{proof}

\begin{lemma}\label{Lemma 3.1}
Let $F=L(H,\chi)$ be a holomorphic line bundle over an abelian surface $A=\BC^2/L$
and let $p:\BC^2\lrt A$ be the natural projection. Then
\vskip 2mm \noindent
{\rm (1)} $T_x^*(F)$ and $F$ have the same Chern classes for each $x\in A$.
\vskip 2mm \noindent
{\rm (2)} For an element $x\in A$, the automorphic factor $J_x$ for the line bundle
$T_x^*(F)$ is given by
\begin{equation*}
  J_x (\alpha,z)=J_{H,\chi}(\alpha,z+a),\quad \alpha\in L,\ z\in\BC^2,
\end{equation*}
where $a$ is an element of $\BC^2$ such that $p(a)=x.$
\end{lemma}
For details, we refer to \cite {MU} or \cite[Proposition 4.14, p.\,127]{Y1}.

\begin{theorem}\label{Theorem 3.1}
Let $A=\BC^2/L$ be an abelian surface with the identity element $e$. Then
there exists a unique holomorphic line bundle $\mathfrak{P}$ over
$A\times \widehat{A}$ which is trivial on $\{ e\}\times \widehat{A}$ and
which satisfies
\begin{equation*}
  \mathfrak{P}|_{A\times \{\xi\}}\cong \mathfrak{P}_{\xi}\qquad
  {\rm for\ all}\ \xi\in \widehat{A},
\end{equation*}
where $\mathfrak{P}_{\xi}$ is the line bundle over $A$ corresponding to
$\xi\in \widehat{A}.$
\end{theorem}
\begin{proof}
See \cite[pp.\,78--80]{MU} and also \cite[pp.\,328--329]{GH}.
\end{proof}

\vskip 2mm
Now we will describe the Poincar{\'e} bundle over $A\times \widehat{A}$
explicitly. For brevity, we put $V=\BC^2.$ We define an Hermitian form $H$
on $V\times \overline{V}^*$ by
\begin{equation*}
  H((z_1,\ell_1),(z_2,\ell_2)):=\ell_1(z_2)+\overline{\ell_2 (z_1)},
\end{equation*}
where $z_1,z_2\in V$ and $\ell_1,\ell_2\in \overline{V}^*.$ We also define the map
$\chi:L\times \widehat{L}\lrt T$ by
\begin{equation*}
  \chi(\alpha,\ell):=\exp\left\{ -i\,\pi\, {\rm Im}\, \ell (\alpha)\right\},
  \quad \alpha\in L,\ \ell\in \widehat{L}.
\end{equation*}
Here $\widehat{L}$ denotes the dual lattice of $L$.
We see easily that $\chi$ satisfies the following equation
\begin{equation*}
  \chi((\alpha+\beta,\ell+\widehat{\ell}\,))=\chi(\alpha,\ell)\,\chi(\beta,\widehat{\ell}\,)
  \,\exp \left\{ i\,\pi\,E((\alpha,\ell),(\beta,\widehat{\ell}\,))\right\},
\end{equation*}
where $\alpha,\beta\in L,\ \ell,\widehat{\ell}\in \widehat{L}$ and $E={\rm Im}\,H.$
Thus $\chi$ is a semi-character of $L\times \widehat{L}$ with respect to $H$.
Then the line bundle $L(H,\chi)$ over $A\times \widehat{A}$ defined by the Hermitian
form $H$ and the semi-character $\chi$ is the Poincar{\'e} bundle over $A\times \widehat{A}$.
In fact, the corresponding automorphic factor
$J_{H,\chi}:(L\times \widehat{L}\,)\times (V\times \overline{V}^*)\lrt GL(1,\BC)$ for
$L(H,\chi)$ is given by
\begin{equation*}
  J_{H,\chi}((\alpha,\widehat{\ell}\,),(z,\ell))= \chi (\alpha,\widehat{\ell}\,)\,
  \exp\,\left\{ \pi\,H((z,\ell),(\alpha,\widehat{\ell}\,))+ \frac{\pi}{2}\,
  H((\alpha,\widehat{\ell}\,)),(\alpha,\widehat{\ell}\,))\right\},
\end{equation*}
where $\alpha\in L,\ \widehat{\ell}\in \widehat{L},\ z\in V$ and $\ell\in \overline{V}^*.$
Let $\widehat{p}:\overline{V}^*\lrt \widehat{A}$ be the natural projection.
The line bundle $L(H,\chi)|_{A\times \{\widehat{p}(\ell)\}}$ over $A$ corresponding to
the point $\widehat{p}(\ell)\in \widehat{A}\,(\ell\in \overline{V}^*)$ is defined by a flat
automorphic factor $J_{\ell}:L\times V\lrt GL(1,\BC)$ given by
\begin{equation*}
  J_{\ell}(\alpha,z)=\exp\,\{\pi\,\ell(\alpha)\}\qquad
  {\rm for\ all}\ \alpha\in L\ {\rm and}\ z\in V.
\end{equation*}
Let $h_{\ell}:V\lrt \BC^*$ ($\ell\in \overline{V}^*$) be a holomorphic function on $V$ defined by
$$h_{\ell}(z):=\exp\,\{-\pi\,\overline{\ell(z)} \},\quad z\in V.$$
Then we have the following relation
\begin{equation*}
  h_\ell (z+\alpha)\,J_\ell (\alpha,z)=\exp\,\{ 2\,\pi\,i\,{\rm Im}\,\ell (\alpha)\}\,h_\ell(z)
\end{equation*}
for all $\alpha\in L$ and $z\in V$. Therefore we obtain the isomorphism
\begin{equation*}
  L(H,\chi)|_{A\times \widehat{p}(\ell)}\cong L(0,\chi_\ell),
\end{equation*}
where $\chi_\ell(\alpha)=\exp\,\{2\,\pi\,i\,\ell(\alpha)\}$ is a semi-character of
$L\,(\alpha\in L)$. We note that if $\ell=0$, then
\begin{equation*}
  L(H,\chi)|_{A\times \{ \widehat{e}\}}\cong L(0,1)\cong \mathcal{O}_A,
\end{equation*}
where $\widehat{e}$ is the identity element of $\widehat{A}$ and
$\mathcal{O}_A$ denotes the trivial line bundle over $A$.
Obviously
\begin{equation*}
  L(H,\chi)|_{\{e\}\times \widehat{A}}\cong L(0,1)\cong \mathcal{O}_{\widehat{A}},
\end{equation*}
where $e$ is the identity element of $A$ and
$\mathcal{O}_{\widehat{A}}$ denotes the trivial line bundle over $\widehat{A}$.
By the uniqueness, $L(H,\chi)$ must be the Poincar{\'e} bundle over $A\times \widehat{A}$.

\vskip 3mm
Now we fix an element $\Om=\begin{pmatrix}
                             \tau & z \\
                             z & \tau
                           \end{pmatrix}\in \widehat{\BH}_2$.
Clearly ${\rm Im}\,\tau> |{\rm Im}\,z|\geq 0.$ Let
$$L_\Om:=\BZ^2\Om+\BZ^2,\quad \BZ^2:=\BZ^{(1,2)}$$
be the lattice in $\BC^2$.
\begin{equation}\label{(3.7)}
  e_1=(\tau,z),\ \,e_2=(z,\tau),\ \,e_3=(1,0)\quad {\rm and}\ e_4=(0,1)
\end{equation}
form a $\BZ$-basis of $L$ in $\BC^2$.
Then the two-dimensional complex torus
\begin{equation*}
  A_\Om:=\BC^2/L_\Om
\end{equation*}
is a {\it special\ abelian\ surface}.

\begin{theorem}\label{Theorem 3.2}
Let $\Om=\begin{pmatrix}
                             \tau & z \\
                             z & \tau
                           \end{pmatrix}\in \widehat{\BH}_2$
such that both ${\rm Im}\,\tau$ and ${\rm Im}\,z$ are integers.
Let $H_\Om:\BC^2\times\BC^2\lrt \BC$ be the Hermitian form on $\BC^2$
defined by
\begin{equation*}
  H_\Om((z_1,z_2),(w_1,w_2)):=z_1\,\overline{w}_1+z_2\,\overline{w}_2,
\end{equation*}
where $(z_1,z_2),(w_1,w_2)\in \BC^2.$
Let $\chi_\Om$ be the semi-character of $L_\Om$ determined uniquely by $H_\Om.$
Then $H_\Om$ is a positive definite Riemann forms for $A_\Om$. The line bundle
$L(H_\Om,\chi_\Om)$ is an ample line bundle over $A_\Om$ and hence the global sections of
$L(H_\Om,\chi_\Om)^{\otimes n}\ (n\geq 3)$ give an embedding of $A_\Om$ into
a projective space. And
\begin{equation*}
\dim H^0(A_\Om,L(H_\Om,\chi_\Om))=\sqrt{\det\,E_\Om}=({\rm Im}\,\tau)^2-({\rm Im}\,z)^2,
\end{equation*}
where $E_\Om$ is the imaginary part of $H_\Om$.
\end{theorem}

\begin{proof}
Let $S_\Om$ (resp.\ $E_\Om$) be the real (resp. imaginary) part of $H_\Om$, that is,
\begin{equation*}
  H_\Om(z,w)=S_\Om(z,w)+i\,E_\Om (z,w)=E_\Om(i\,z,w)+i\,E_\Om(z,w),\qquad z,w\in \BC^2.
\end{equation*}
We note that
$$
S_\Om (z,w)=S_\Om(w,z)=E_\Om(i\,z,w)\ \,{\rm and}\ \,E_\Om(z,w)=-E_\Om(w,z),
\quad z,w\in\BC^2.
$$
Let $e_1,\,e_2,\,e_3$ and $e_4$ be the basis of $L_\Om$ given by
\eqref{(3.7)}. The $4\times 4$
matrices $\mathbb{S}_\Om:=(S_\Om(e_i,e_j))_{1\leq i,j\leq 4}$ and
$\mathbb{E}_\Om:=(E_\Om(e_i,e_j))_{1\leq i,j\leq 4}$ are given by
\begin{equation*}
\mathbb{S}_\Om=
\begin{pmatrix}
  |\tau|^2+|z|^2 & \tau\, \overline{z}+\overline{\tau}\,z & {\rm Re}\,\tau & {\rm Re}\,z \\
  \tau\, \overline{z}+\overline{\tau}\,z & |\tau|^2+|z|^2 & {\rm Re}\,z & {\rm Re}\,\tau \\
  {\rm Re}\,\tau & {\rm Re}\,z & 1 & 0 \\
  {\rm Re}\,z & {\rm Re}\,\tau & 0 & 1
\end{pmatrix}
\end{equation*}
and
\begin{equation*}
\mathbb{E}_\Om=
  \begin{pmatrix}
 \ 0 & \ 0 & {\rm Im}\,\tau & {\rm Im}\,z \\
 \ 0 & \ 0 & {\rm Im}\,z & {\rm Im}\,\tau \\
  -{\rm Im}\,\tau & -{\rm Im}\,z & 1 & 0 \\
  -{\rm Im}\,z & -{\rm Im}\,\tau & 0 & 1
\end{pmatrix}.
\end{equation*}
Since both ${\rm Im}\,\tau$ and ${\rm Im}\,z$ are integers, $\mathbb{E}_\Om$ is an integral
alternating matrix. Therefore $E_\Om$ is integral valued on $L_\Om\times L_\Om.$
It is obvious that $H_\Om$ is positive definite. Therefore $H_\Om$ is a positive definite
Riemann form for $A_\Om$ and hence $L(H_\Om,\chi_\Om)$ is an ample line bundle over
$A_\Om$. According to Theorem of Lefschetz\,(cf.\,\cite[pp.\,29--30]{MU}),
the global sections of
$L(H_\Om,\chi_\Om)^{\otimes n}\ (n\geq 3)$ gives an embedding of $A_\Om$ into
a projective space. By an easy calculation, we get
\begin{equation*}
  \det \mathbb{E}_\Om=\left\{ ({\rm Im}\,\tau)^2-({\rm Im}\,z)^2 \right\}^2.
\end{equation*}
According to Formula \eqref{(3.5)}, we get
\begin{equation*}
\dim H^0(A_\Om,L(H_\Om,\chi_\Om))=\sqrt{\det\,E_\Om}=({\rm Im}\,\tau)^2-({\rm Im}\,z)^2.
\end{equation*}
\end{proof}

\begin{theorem}\label{thm:3.3}
Let $\Om=\begin{pmatrix}
                             \tau & 0 \\
                             0 & \tau
                           \end{pmatrix}\in \widehat{\BH}_2$.
Let $H_\tau:\BC^2\times\BC^2\lrt \BC$ be the Hermitian form on $\BC^2$
defined by
\begin{equation*}
  H_\tau((z_1,z_2),(w_1,w_2)):=\frac{1}{{\rm Im}\,\tau}
  (z_1\,\overline{w}_1+z_2\,\overline{w}_2),\qquad
 {\rm where}\ (z_1,z_2),(w_1,w_2)\in \BC^2.
\end{equation*}
Let $\chi_\tau$ be the semi-character of $L_\Om$ with respect to $E_\tau={\rm Im}\,H_\tau$.
Then $H_\tau$ is a positive definite Riemann forms for $A_\Om$. The line bundle
$L(H_\tau,\chi_\tau)$ is an ample line bundle over $A_\Om$ and hence the global sections of
$L(H_\tau,\chi_\tau)^{\otimes n}\ (n\geq 3)$ gives an embedding of $A_\Om$ into
a projective space. And
\begin{equation*}
\dim H^0(A_\Om,L(H_\tau,\chi_\tau))=1.
\end{equation*}
\end{theorem}

\begin{proof}
Let $S_\tau$ (resp.\ $E_\tau$) be the real (resp. imaginary) part of $H_\tau$, that is,
\begin{equation*}
  H_\tau(z,w)=S_\tau(z,w)+i\,E_\tau (z,w)=E_\tau(i\,z,w)+i\,E_\tau(z,w),\qquad z,w\in \BC^2.
\end{equation*}
We note that
$$
S_\tau (z,w)=S_\tau(w,z)=E_\tau(i\,z,w)\ \,{\rm and}\ \,E_\tau(z,w)=-E_\tau(w,z),
\quad z,w\in\BC^2.
$$
Let
$$f_1=(\tau,0),\,f_2=(0,\tau),\,f_3=(1,0)\ {\rm and}\ f_4=(0,1)$$
be a $\BZ$-basis of $L_\Om$ .
The $4\times 4$ matrices $\mathbb{S}_\tau:=(S_\tau(f_i,f_j))_{1\leq i,j\leq 4}$ and
$\mathbb{E}_\tau:=(E_\Om(f_i,f_j))_{1\leq i,j\leq 4}$ are given by
\begin{equation*}
\mathbb{S}_\tau=\frac{1}{{\rm Im}\,\tau}\,
\begin{pmatrix}
  |\tau|^2 & 0 & {\rm Re}\,\tau & 0 \\
  0 & |\tau|^2 & 0 & {\rm Re}\,\tau \\
  {\rm Re}\,\tau & 0 & 1 & 0 \\
  0 & {\rm Re}\,\tau & 0 & 1
\end{pmatrix}
\end{equation*}
and
\begin{equation*}
\mathbb{E}_\tau=\frac{1}{{\rm Im}\,\tau}\,
  \begin{pmatrix}
 \ 0 & \ 0 & {\rm Im}\,\tau & 0 \\
 \ 0 & \ 0 & 0 & {\rm Im}\,\tau \\
  -{\rm Im}\,\tau & 0 & 0 & 0 \\
  \ 0 & -{\rm Im}\,\tau & 0 & 0
\end{pmatrix}
= \begin{pmatrix}
 \ 0 & \ 0 & 1 & 0 \\
 \ 0 & \ 0 & 0 & 1 \\
  -1 & \ 0 & 0 & 0 \\
  \ 0 & -1 & 0 & 0
\end{pmatrix}.
\end{equation*}
Clearly $\mathbb{E}_\tau$ is an integral alternating matrix and so $E_\tau$ is
integral valued on $L_\Om\times L_\Om.$
It is obvious that $H_\tau$ is positive definite. Therefore $H_\tau$ is a positive definite
Riemann form for $A_\Om$ and hence $L(H_\tau,\chi_\tau)$ is an ample line bundle over
$A_\Om$. According to Theorem of Lefschetz\,(cf.\,\cite[pp.\,29--30]{MU}),
the global sections of $L(H_\tau,\chi_\tau)^{\otimes n}\ (n\geq 3)$ gives an embedding
of $A_\Om$ into a projective space. According to Formula \eqref{(3.5)}, we get
\begin{equation*}
\dim H^0(A_\Om,L(H_\tau,\chi_\tau))=\sqrt{\det\,E_\tau}=1.
\end{equation*}
\end{proof}

\begin{theorem}\label{thm:3.4}
Let $\Om=\begin{pmatrix}
                             \tau & z \\
                             z & \tau
                           \end{pmatrix}\in \widehat{\BH}_2$
such that both ${\rm Im}\,(2\,\tau+z)$ and ${\rm Im}\,(\tau+2\,z)$ are integers.
Let $H_*:\BC^2\times\BC^2\lrt \BC$ be the Hermitian form on $\BC^2$
defined by
\begin{equation*}
  H_*((z_1,z_2),(w_1,w_2)):=2\,z_1\,\overline{w}_1+2\,z_2\,\overline{w}_2
  +z_1\,\overline{w}_2+z_2\,\overline{w}_1,
\end{equation*}
where $(z_1,z_2),(w_1,w_2)\in \BC^2$.
Then $H_*$ is a positive definite Riemann forms for $A_\Om$.
Let $\chi_*$ be the semi-character of $L_\Om$ determined uniquely by $H_*$. Then
the line bundle $L(H_*,\chi_*)$ is an ample line bundle over $A_\Om$ and
hence the global sections of $L(H_*,\chi_*)^{\otimes n}\ (n\geq 3)$ gives
an embedding of $A_\Om$ into a projective space. And
\begin{equation*}
\dim H^0(A_\Om,L(H_*,\chi_*))=\sqrt{\det\,E_*}=
3\,\left\{ ({\rm Im}\,\tau)^2-({\rm Im}\,z)^2\right\},
\end{equation*}
where $E_*$ is the imaginary part of $H_*$.
\end{theorem}

\begin{proof}
Let $S_*$ (resp.\ $E_*$) be the real (resp. imaginary) part of $H_*$, that is,
\begin{equation*}
  H_*(z,w)=S_*(z,w)+i\,E_* (z,w)=E_*(i\,z,w)+i\,E_*(z,w),\qquad z,w\in \BC^2.
\end{equation*}
We note that
$$
S_* (z,w)=S_*(w,z)=E_*(i\,z,w)\ \,{\rm and}\ \,E_*(z,w)=-E_*(w,z),
\quad z,w\in\BC^2.
$$
Let $e_1,\,e_2,\,e_3$ and $e_4$ be the $\BZ$-basis of $L_\Om$ given by \eqref{(3.7)}.
The $4\times 4$ matrices $\mathbb{S}_*:=(S_*(e_i,e_j))$ and
$\mathbb{E}_*:=(E_*(e_i,e_j))$ are given by
\begin{equation*}
\mathbb{S}_*=
\begin{pmatrix}
  2\,(|\tau|^2+|z|^2) +\tau\,\overline{z}+z\,\overline{\tau}
  & |\tau|^2+|z|^2+2\,(\tau\, \overline{z}+\overline{\tau}\,z)
  & {\rm Re}\,(2\,\tau+z) & {\rm Re}\,(\tau+2\,z) \\
  |\tau|^2+|z|^2+2\,(\tau\, \overline{z}+\overline{\tau}\,z)  &
  2\,(|\tau|^2+|z|^2) +\tau\,\overline{z}+z\,\overline{\tau} &
  {\rm Re}\,(\tau+2\,z) & {\rm Re}\,(2\,\tau+z) \\
  {\rm Re}\,(2\,\tau+z) & {\rm Re}\,(\tau+2\,z) & 2 & 1 \\
  {\rm Re}\,(\tau+2\,z) & {\rm Re}\,(2\,\tau+z)& 1 & 2
\end{pmatrix}
\end{equation*}
and
\begin{equation*}\
\mathbb{E}_*=
  \begin{pmatrix}
 \ 0 & \ 0 & {\rm Im}\,(2\,\tau+z) & {\rm Im}\,(\tau+2\,z) \\
 \ 0 & \ 0 & {\rm Im}\,(\tau+2\,z) & {\rm Im}\,(2\,\tau+z) \\
  -{\rm Im}\,(2\,\tau+z) & -{\rm Im}\,(\tau+2\,z) & 0 & 0 \\
  -{\rm Im}\,(\tau+2\,z) & -{\rm Im}\,(2\,\tau+z) & 0 & 0
\end{pmatrix}.
\end{equation*}
Since both ${\rm Im}\,(2\,\tau+z)$ and ${\rm Im}\,(\tau+2\,z)$ are integers,
$\mathbb{E}_*$ is an integral alternating matrix. We see that $E_*$ is
integral valued on $L_\Om\times L_\Om.$
It is obvious that $H_*$ is positive definite. Therefore $H_*$ is a positive definite
Riemann form for $A_\Om$ and hence $L(H_*,\chi_*)$ is an ample line bundle over
$A_\Om$. According to Theorem of Lefschetz\,(cf.\,\cite[pp.\,29--30]{MU}),
the global sections of $L(H_*,\chi_*)^{\otimes n}\ (n\geq 3)$ gives an embedding
of $A_\Om$ into a projective space. According to Formula \eqref{(3.5)},
\begin{equation*}
\dim H^0(A_\Om,L(H_*,\chi_*))=\sqrt{\det\,E_*}=
3\,\left\{ ({\rm Im}\,\tau)^2-({\rm Im}\,z)^2 \right\}.
\end{equation*}
\end{proof}

\vskip 3mm
Let
$\Om=\begin{pmatrix}
\tau & z \\
z & \tau
\end{pmatrix}\in \widehat{\BH}_2$ and let $L_\Om=\BZ^2\Om+\BZ^2$ be the
lattice in $\BC^2$. Then $A_\Om=\BC^2/L_\Om$ is a special abelian surface.
We choose a positive definite Riemann form $H_\Om$ for $A_\Om$ with the
semi-character $\chi_\Om$ of the lattice $L_\Om$ determined by $H_\Om$.
We have the automorphic factor $J_\Om:L_\Om\times \BC^2 \lrt GL(1,\BC)$
given by
\begin{equation*}
  J_\Om(\alpha,z)=\chi_\Om (\alpha)\,\exp\,\left\{\pi\,H_\Om(z,\alpha)+
  \frac{\pi}{2}\,H_\Om(\al,\al) \right\},\qquad \al\in L_\Om,\ z\in\BC^2.
\end{equation*}
See the formula \eqref{(3.2)}. For each $\al\in L_\Om$, we define the holomorphic function
$J_{\Om,\al}:\BC^2\lrt \BC$ by
\begin{equation}\label{(3.8)}
  J_{\Om,\al}(z):=J_\Om (\al,z),\quad z\in\BC^2.
\end{equation}
Then we can show that the mapping
\begin{equation}\label{(3.9)}
 L_\Om \lrt H^0 (\BC^2,\mathcal{O}_{\BC^2}^*),\quad
 \al\longmapsto J_{\Om,\al}\ (\al\in L_\Om)
\end{equation}
is a 1-cocycle on $L_\Om$ with coefficients in $H^0 (\BC^2,\mathcal{O}_{\BC^2}^*)$
defining the Chern class of $L(H_\Om,\chi_\Om)$ being $E_\Om\in H^2(A_\Om,\BZ).$
Here $\mathcal{O}_{\BC^2}$ denotes the structure sheaf on $\BC^2$ and $E_\Om$
denotes the imaginary part of $H_\Om$.
\vskip 3mm
We review the geometric properties of $A_\Om$. For an integer $r$ with $0\leq r\leq 4$,
we denote by $\mathfrak{A}_r(\Om)$ the group of all alternating $r$-forms
$L_\Om\times\cdots\times L_\Om\lrt \BZ$. Let $\Om^q(A_\Om)$ denote the sheaf of
holomorphic $q$-forms on $A_\Om$. Then we have the following properties\,:
\vskip 2mm\noindent
($\Om$1) {\it We have canonical isomorphisms}
\begin{equation}\label{(3.10)}
  H^r(A_\Om,\BZ)\cong \mathfrak{A}_r(\Om),\quad r=0,1,2,3,4.
\end{equation}

\vskip 2mm\noindent
($\Om$2) {\it Let $V$ be the tangent space to $A_\Om$ at $e$.
Here $e$ is the identity element $e$ of $A_\Om.$
Let $T={\rm Hom}_\BC (V,\BC)$ be the complex cotangent space to
$A_\Om$ at $e$. We put
$$\overline{T}={\rm Hom}_{\BC\!-\!{\rm antilinear}}(V,\BC).$$
Then there are natural isomorphisms
\begin{equation*}
  H^q(A_\Om,\mathcal{O}_\Om)\cong \Lambda^q \overline{T}
\end{equation*}
for all $q$ and hence
\begin{equation*}
  H^q(A_\Om,\Om^p(A_\Om))\cong \Lambda^p T\otimes \Lambda^q \overline{T}.
\end{equation*}
Here $\mathcal{O}_\Om$ denotes the structure sheaf on $A_\Om$.}
\vskip 2mm
We refer to \cite[pp.\,3--8]{MU} for the proofs of ($\Om$1) and ($\Om$2)

\vskip 3mm
From the following exact sequence
$$ 0\lrt \BZ \lrt \mathcal{O}_\Om
\stackrel{e^{2\pi i(\cdot)}}{\lrt} \mathcal{O}_\Om^* \lrt 0,$$
we get the exact sequence
\begin{equation*}
  0\lrt H^0(A_\Om,\BZ) \lrt  \cdots
  \lrt H^1(A_\Om,\mathcal{O}_\Om^*) \stackrel{c_1}{\lrt} H^2(A_\Om,\BZ)
  \stackrel{\mathfrak{K}}{\lrt}
  H^2(A_\Om,\mathcal{O}_\Om)\lrt \cdots.
\end{equation*}
By ($\Om$1), we have $H^2(A_\Om,\BZ)\cong \mathfrak{A}_2(\Om).$
The alternation 2-form $E_\Om$ is nothing but the image of the first cohomology class
defined by the 1-cocycle \eqref{(3.9)} under the mapping $c_1$. Therefore $E_\Om$ is the
first Chern class of $L(H_\Om,\chi_\Om)$, that is,
$E_\Om=c_1(L(H_\Om,\chi_\Om))\in H^2(A_\Om,\BZ)$ (cf. Formula \eqref{(3.10)}).
The Picard group ${\rm Pic}\,(A_\Om)$ of $A_\Om$ is defined to be the group of
isomorphism classes of line bundles over $A_\Om$ with group operation being tensor
product. Then we have the following exact sequence
\begin{equation*}
  1\lrt {\rm Pic}^0\,(A_\Om)\lrt {\rm Pic}\,(A_\Om)\lrt {\rm Ker}\,\mathfrak{K}
  \lrt 1,
\end{equation*}
where ${\rm Ker}\,\mathfrak{K}$ denotes the kernel of the mapping
$\mathfrak{K}: H^2(A_\Om,\BZ)\lrt H^2(A_\Om,\mathcal{O}_\Om)$.
Since
\begin{equation*}
  \widehat{A}_\Om={\rm Pic}^0\,(A_\Om)\cong {\rm Hom}\,(L_\Om,T)
  \quad {\rm and}\quad {\rm Pic}\,(A_\Om) \cong H^1(A_\Om,\mathcal{O}_\Om^*),
\end{equation*}
we have the exact sequence
\begin{equation*}
1\lrt {\rm Hom}\,(L_\Om,T)(=\widehat{A}_\Om)\lrt H^1(A_\Om,\mathcal{O}_\Om^*)\lrt
{\rm Ker}\,\mathfrak{K}\lrt 1.
\end{equation*}
We recall that $\widehat{A}_\Om$ is the dual abelian surface of $A_\Om$.

\vskip 3mm
We now compute the curvature of the line bundle $\mathfrak{L}_\Om=L(H_\Om,\chi_\Om)$
over $A_\Om$. We fix an Hermitian structure $h_\Om$ on $\mathfrak{L}_\Om.$
We pull back $h_\Om$  to $p^*\mathfrak{L}_\Om=\widetilde{\mathfrak{L}}_\Om$
to obtain an Hermitian structure $\widetilde{h}_\Om$
on $\widetilde{\mathfrak{L}}_\Om=\BC^2\times\BC.$
We recall that $p:\BC^2\lrt A_\Om$ is the natural projection.
We may consider $\widetilde{h}_\Om$ as a positive function on $\BC^2$ invariant
under the action of $L_\Om$. That is, $\widetilde{h}_\Om$ satisfies
\begin{equation*}
  \widetilde{h}_\Om(z)=|J_{\al,z}|^2\,\widetilde{h}_\Om (z+\al),\quad
  \al\in L_\Om,\ z\in \BC^2.
\end{equation*}
Then the connection form $\widetilde{\theta}(z)=\partial\log\,\widetilde{h}_\Om$
and the curvature form
$\widetilde{\Theta}=\overline{\partial}\partial\log\,\widetilde{h}_\Om$ are given by
\begin{equation*}
  \widetilde{\theta}(z)=\widetilde{\theta}(z+\al)+\partial \log J_\al (\al,z)
  \quad {\rm and}\quad \widetilde{\Theta}(z)=\widetilde{\Theta}(z+\al).
\end{equation*}
Therefore the curvature form of $\mathfrak{L}_\Om$ is an ordinary 2-form on $A_\Om$.
Multiplying $H_\Om$ by a suitable positive $C^\infty$-function on $A_\Om$, we may assume
that the curvature form of $\mathfrak{L}_\Om$ is a harmonic (1,1)-form on $A_\Om$.
Since
$$H^2(A_\Om,\BC)=\Lambda^2 V^*\otimes V^*\oplus \overline{V}^*\oplus
\Lambda^2 \overline{V}^*,\quad V:=\BC^2,$$
a harmonic form on $A_\Om$ has constant coefficients with respect to  the natural
coordinate $z_1,z_2$ in $V$. In particular, the curvature $\Theta$ of
$\mathfrak{L}_\Om=L(H_\Om,\chi_\Om)$ is given by
\begin{equation*}
  \Theta=\sum_{i,j=1}^2\,c_{ij}\,dz_i\wedge d\overline{z}_j,
\end{equation*}
where the $c_{ij}$ are constant.

\end{section}

\vskip 10mm

\end{document}